\documentclass[a4paper,12pt]{article}
\usepackage{cmap}                        
\usepackage[cp1251]{inputenc}            
\usepackage[english]{babel}
\usepackage[left=2cm,right=18mm,top=15mm,bottom=20mm]{geometry} 
\usepackage{amsmath, amssymb}

\begin{document}

\begin{center} \textbf{ ONE FORMATION  OF FINITE GROUPS}

\vspace{+1mm}

\emph{Viachaslau I. Murashka}
\footnote{Francisk Skorina Gomel State University,
Gomel, Belarus,
e-mail: mvimath@yandex.ru}\end{center}

\vspace{-2mm}

\textbf{Abstract.} It is shown that two formations of finite groups, one was introduced
 by V.S. Monakhov and V.N. Kniahina and another one was introduced by R. Brandl, are coincides.

\smallskip

In this note only finite groups are considered. V.S. Monakhov and V.N. Kniahina \cite{3}
 studied the class
 $\mathfrak{X}$ of all groups whose cyclic primary subgroups are  $\mathbb{P}$-subnormal.
  It was shown that  $\mathfrak{X}$  is a hereditary saturated formation and
 $\mathfrak{U}\subset \mathfrak{X}\subset\mathfrak{D}$ where $\mathfrak{U}$
  and $\mathfrak{D}$ are classes of all supersoluble groups and groups with Sylow tower
 of supersoluble type respectively. In \cite{6, 7} R. Brandl studied groups satisfying some
 law $\ddot{u}_k(x,y)=1$ where
 $\ddot{u}_1(x,y)=[x,y]$ and
  $\ddot{u}_{k+1}(x,y)=\ddot{u}_{k}(x,y)^{-k}[\ddot{u}_{k}(x,y),y]$ for $k>1$.
 He showed that the class $\mathfrak{B}=(G|$ for all $x,y\in G$ there is a natural
 $k$ such that $\ddot{u}_k(x,y)=1)$ is a hereditary saturated formation containing
  $\mathfrak{U}$. Also he showed that this class coincides with the class of all
 groups whose subgroups with nilpotent derived subgroup are supersoluble. In this note
 will be proved that  classes  $\mathfrak{B}$ and $\mathfrak{X}$ are coincides.

 The standard definitions and notation from \cite{5} are used. Recall  \cite{1} that a subgroup $H$ of
 a group $G$ is called $\mathbb{P}$-subnormal if either  $H=G$ or there is a maximal chain of subgroups  $H=H_0\subset H_1\subset\dots\subset H_n = G$ such that $|H_i:H_{i-1}|$ is a prime for  $i = 1,\dots,n$.

\textbf{Theorem.} \emph{For a group $G$ the following statements are equivalent}:

(1) \emph{All cyclic primary subgroups of $G$ are $\mathbb{P}$-subnormal in $G$}.

(2) \emph{All subgroups of $G$ with nilpotent derived subgroup are supersoluble}.

(3) \emph{For all $x,y\in G$ there is natural $k$ such that  $\ddot{u}_k(x,y)=1$}.

\emph{Proof.} The equivalence of Assertions (2) and (3) was established in $\cite{6}$. Also there
 were shown that $\mathfrak{B}=LF(f)$ where  $f(p)$ is the class of all soluble groups of
 exponent dividing  $p-1$ for all primes $p$. It is well known that $f(p)$ is hereditary
 formation for all primes $p$.

Assume that $\mathfrak{X}$ is not contained in $\mathfrak{B}$. Let choose a group  $G$
 of minimal order from  $\mathfrak{X}\setminus \mathfrak{B}$. Note \cite{1} that $G$ is soluble.
 Since $\mathfrak{X}$ and $\mathfrak{B}$ are saturated formations, we see that  $\Phi(G) = 1$.
 It is clear that  $G$ is $\mathfrak{B}$-critical group. Since $\mathfrak{B}$ is a formation,
 we note that $G$ has unique minimal normal subgroup  $N$. Now $N$  is abelian $p$-subgroup
 and $N=C_G(N)$. Let $ M$ be a maximal subgroup $G$ which does not contain $N$. Then $G = NM$.
 Since $N$ is abelian,   $N\cap M = 1$. Since $G$ has Sylow tower of supersoluble type, $p$ is
 the maximal prime divisor of $|G|$ and $(p,|M|) = 1$. Let $H$ be a proper subgroup of  $M$.
  From  $N=C_G(N)$ it follows that $O_{p'}(NH)=1$. Now $O_{p',p}(NH)=N$. From $NH\in\mathfrak{B}$
 it follows that $H\simeq NH/O_{p',p}(NH)\in f(p)$. Hence $M$ is $f(p)$-critical group. If
 $M$ is not a cyclic primary subgroup then  $M\in f(p)$, a contradiction. Thus $M$ is a cyclic
 primary subgroup. By theorem B of \cite{3} $G$ is supersoluble, a contradiction. So
 $\mathfrak{X}\subseteq \mathfrak{B}$.

Assume now that $\mathfrak{B}$ is not contained in $\mathfrak{X}$. Let choose a group $G$
 of minimal order from  $\mathfrak{B}\setminus \mathfrak{X}$. It is easy to show that
 $\Phi(G) = 1$ and $G$ is an $\mathfrak{X}$-critical group. By \cite{3} $G$ is a biprimary minimal
 non-supersoluble group with unique minimal normal subgroup $N$ which is the Sylow
  $p$-subgroup, a Sylow $q$-subgroup $Q$ of $G$ is cyclic. It means that the derived subgroup of $G$
 is nilpotent. Hence $G$ is supersoluble. Now $G\in\mathfrak{X}$, the final contradiction.
 So $\mathfrak{B}\subseteq\mathfrak{X}$. Thus $\mathfrak{B} = \mathfrak{X}$. $\square$


\vspace{-7mm}


\begin{thebibliography}{99}\small\vspace{-3mm}


\bibitem{3}
 V.S. Monakhov and V.N. Kniahina, Finite groups with $\mathbb{P}$-subnormal subgroups,
 Ricerche di Matematica, \textbf{DOI} 10.1007/s11587-013-0153-9, 2013.\vspace{-4mm}

 \bibitem{6}
R. Brandl, Groups sharing some varietal properties with supersoluble groups, J. Austral.
 Math. Soc. \textbf{34}
  (1981), 265-268.\vspace{-4mm}

\bibitem{7}
 R. Brandl,
  Zur Theorie der untergruppenabgeschlossenen Formationen: Endliche Variet\"{a}ten, J.
 Algebra,\textbf{ 73}
  (1981),  1-22.\vspace{-4mm}



\bibitem{5}
 K. Doerk and T. Hawkes, Finite soluble groups,   Walter de Gruyter,
 1992.\vspace{-4mm}



\bibitem{1}   A.F. Vasil’ev,
T.I. Vasil’eva and
V.N. Tyutyanov,  On the finite groups of supersoluble type, Siberian Mathematical
 Journal,  \textbf{51}(6) (2010),  1004-1012.
\end{thebibliography}
\end{document}